# Definitions of real order integrals and derivatives using operator approach

## Raoelina Andriambololona

Theoretical Physics Dept., Antananarivo, Madagascar, Institut National des Sciences et Techniques Nucléaires (INSTN-Madagascar), Boite Postale 4279, Antananarivo 101, Madagascar

**Email address:**
raoelinasp@yahoo.fr (R. Andriambololona)

**Abstract:** The set $E$ of functions $f$ fulfilling some conditions is taken to be the definition domain of $s$-order integral operator $J^s$ (iterative integral), first for any positive integer $s$ and after for any positive $s$ (fractional, transcendental $\pi$ and $e$). The definition of $k$-order derivative operator $D^k$ for any positive $k$ (fractional, transcendental $\pi$ and $e$) is derived from the definition of $J^s$. Some properties of $J^s$ and $D^k$ are given and demonstrated. The method is based on the properties of Euler's gamma and beta functions.

## 1. Introduction

The definition of the first order ordinary integral (respectively derivative) of an integrable (respectively derivable) function $f$ of a real variable $t$ is the existence of

$$J^1(f)(x) = \lim_{|t_{n+1}-t_n|\to 0} \sum_{t_n=a}^{x} f(t_n)(t_{n+1}-t_n)$$

noted

$$\int_a^x f(t)dt$$

$$D^1(f)(x) = \lim_{|t-x|\to 0} \frac{f(t)-f(x)}{t-x}$$

written $\frac{df}{dx}$ in Leibniz's notation and $\frac{d}{dx}(f)(x)$ in operator notation. By mathematical induction, we extend the definition of $J^s(f)(x)$ and $D^s(f)(x)$ for any positive integer $s$.

The problem of fractional integrals and derivatives is: "may we extend the definitions of $J^s$ and $D^s$ for $s$ positive fraction, negative fraction, real and complex numbers ?".

In fact, the fractional derivatives problem is an old one. In 1695, Leibniz raised the question: "Can the meaning of derivatives with integer order be generalized to derivatives with non-integer orders? "L'Hospital replied by another question to Leibniz:" What if the order be 1/2 ? "And Leibniz replied: "It will lead to a paradox, from which one day useful consequences will be drawn." This was the birth of fractional derivatives.

Since then, several approaches have been done [1],[2],[3],[4] In reference [5], we have considered the case of $f(x) = ax^k$ and for $x \in \mathbb{R}_+$. A first definition of a $s$-order derivative is given by the fundamental relation

$$\frac{d^s}{dx^s}(ax^k) = a\frac{\Gamma(k+1)}{\Gamma(k-s+1)}x^{k-s} + \sum_{m=-\infty}^{-1} c_m \frac{x^{m-s}}{\Gamma(m-s+1)}$$

in which $\{c_m\}$ is a set of arbitrary constants with finite values and $\Gamma$ are Euler's gamma function. The second part of this relation is the terms which allow the unified expression of derivatives and integrals respectively for both positive and negative orders. If $s$ is a positive integer, the expression gives the ordinary derivative. If $s$ a negative integer number, the expression is the ordinary indefinite integral (primitive) of $ax^k$; it is not unique because it depends on arbitrary constants. This result is expected because it is well known that an indefinite integral is not unique. If the order $s$ is a negative fraction, the fractional derivative is not unique; it depends on arbitrary constants. Notions of linear, semi-linear, commutative and semi-commutative properties fractional derivatives are introduced too.

A second definition of $s$-order derivative is given. It gives exactly the same results as in the first definition in the case of integer order. But some differences appear between the two definitions in the case of fractional order.

In the present work, a set $E$ of more general functions is used instead of $ax^k$. $E$ is the set of functions $f$ of real variable $x$ such as $f(x) = 0$ for $x \leq 0$. $f$ is derivable for any order and integrable for any order in the interval $]0, +\infty[$. The $s$-order integral operator $J^s$ for $s$ positive integer is defined by $s$-iterative integral of $J^1$. The definition of $J^s$ for any positive $s$ is obtained in extending the expression of $J^s$ for $s$ positive integer to any $s$ positive number. The definition of $k$-order derivative operator $D^k$ for any positive $k$ is derived from the expression of $J^s$ for positive $s$. Properties of $J^s$ and $D^k$ for any positive $s$ and positive $k$ are studied.

New and remarkable results for $J^\pi, J^e, D^\pi$ and $D^e$ for transcendental numbers $\pi$ and $e$ are given too.

In all our work, we define the extension of the factorial function $n!$ as

$$n! = \Gamma(n+1) \text{ for any } n \in \mathbb{R} - \{k / k \in \mathbb{Z}_-\}$$

For instance

$$\left(-\frac{1}{2}\right)! = \Gamma\left(-\frac{1}{2} + 1\right) = \Gamma\left(\frac{1}{2}\right) = \frac{\sqrt{\pi}}{2}$$

$$(e - \pi)! = \Gamma(e - \pi + 1)$$

## 2. Definition of One-Order Operator $J^1$ and One-Order Derivative Operator $D^1$

Let $E$ be the set of function $f$ of real variable $x$, derivable and integrable infinitely (at any order), and verifying

$$f^{(k)}(x) = 0 \text{ for } x \leq 0 \tag{2.1}$$

for any $k \in \mathbb{N}$. $f^{(k)}$ stands for $k$-order derivative of $f$.
Let us define the one order integral operator $J^1$ and the one-order derivative operator $D^1$ by the relations

$$J^1(f)(x) = \int_0^x f(t)\, dt \tag{2.2}$$

$$D^1(f)(x) = \frac{d}{dx}(f)(x) = f'(x) \tag{2.3}$$

$f'(x)$ is the first derivative of $f(x)$. $D^1$ is obviously an operator over $E$ because there is one and only one $D^1(f)(x)$ for a given $f$. $J^1$ is an operator over $E$ too because we have a definite integral. It will not be the case if we have an indefinite integral: $J^1(f)(x)$ will not be unique because it will depend on an arbitrary additive constant.

## 3. Relations between $J^1$ and $D^1$

Theorem 1

$$J^1 D^1 = D^1 J^1 = 1_E \tag{3.1}$$

in which $1_E$ is the identity operator over the set $E$.
Proof

$$J^1 D^1(f)(x) = \int_0^x f'(t)\, dt = f(x) - f(0) = f(x) \tag{3.2}$$

$$J^1 D^1 = 1_E \tag{3.3}$$

$$D^1 J^1(f)(x) = \frac{d}{dx} \int_0^x f(t)\, dt = f(x) \tag{3.4}$$

$$D^1 J^1 = 1_E \tag{3.5}$$

Theorem 2
$J^1$ and $D^1$ are inverse of each other: $D^1$ is the inverse of $J^1$ and $J^1$ is the inverse of $D^1$.
Proof
It is obvious because $J^1 D^1 = D^1 J^1 = 1_E$
Remarks.
We have $D^1 J^1 = 1_E$ but we do not have $J^1 D^1 = 1_E$ if and only if $f(0) \neq 0$.

$$D^1 J^1(f)(x) = \frac{d}{dx} \int_0^x f(t)\, dt = \frac{d}{dx}[F(x) - F(0)]$$

$$= \frac{d}{dx}[F(x)] = f(x) \tag{3.6}$$

$F(x)$ is a primitive of $f(x)$.

$$J^1 D^1(f)(x) = \int_0^x f'(t)\, dt = f(x) - f(0)$$

$$\neq f(x) \text{ if } f(0) \neq 0 \tag{3.7}$$

Example.
For instance, let us take $f(x) = e^x$:

$$D^1 J^1(e^x) = e^x \quad J^1 D^1(e^x) = e^x - 1 \neq e^x$$

1) In general case $J^1$ is the right hand side inverse $D_R^{-1}$ of $D^1$ and $D^1$ is the left hand side inverse $J_L^{-1}$ of $J^1$. (see annex 2)[6]
2) For trigonometric function, the derivative operator is a rotation of $\frac{\pi}{2}$ angle and integral operator is a rotation of angle $-\frac{\pi}{2}$ :

$$D^1(\sin)(x) = \sin\left(x + \frac{\pi}{2}\right) \quad J^1(\sin)(x) = \sin\left(x - \frac{\pi}{2}\right)$$

$$D^s(\sin)(x) = \sin\left(x + s\frac{\pi}{2}\right) \quad J^s(\sin)(x) = \sin\left(x - s\frac{\pi}{2}\right)$$

$$D^1(\cos)(x) = \cos\left(x + \frac{\pi}{2}\right) \quad J^1(\cos)(x) = \cos\left(x - \frac{\pi}{2}\right)$$

$$D^s(\cos)(x) = \sin\left(x + s\frac{\pi}{2}\right) \quad J^s(\cos)(x) = \sin\left(x - s\frac{\pi}{2}\right)$$

for any positive integer $s$. It may be extended for any real and complex $s$.

This property will give us the possibility to extend the definition of the $s$-order operator derivative $D^s$ and the $s$-order integral operator defined for positive integer $s$ to negative integer $s$.

$$J^{-1} = D^1 \text{ is the inverse of } J^1$$
$$D^{-1} = J^1 \text{ is the inverse of } D^1$$
$$J^{-s} = D^s \text{ for any positive } s$$
$$D^{-s} = J^s \text{ for any positive } s$$

Theorem 3
$$D^s J^s = J^s D^s = 1_E \text{ for any positive } s \quad (3.8)$$

Proof
It may be deduced easily from the theorem 1.

## 4. Extension 1. Definition of S-Order Integral Operator $J^s$ and S-Order Derivative Operator $D^s$ for any Positive Integer $s$.

Let us iterate $s$-times $J^1$

$$J^0(f)(x) = f(x)$$

$$J^1(f)(x) = \int_0^x f(t)dt$$

$$J^2(f)(x) = \int_0^x J^1(f)(t)dt = \int_0^x \left[\int_0^{t_1} (f)(t_2)dt_2\right] dt_1$$

$$\ldots\ldots\ldots\ldots\ldots$$

$$J^s(f)(x) = \int_0^x \int_0^{t_1} \int_0^{t_2} \ldots \int_0^{t_{s-1}} f(t_s)\, dt_s \ldots dt_2 dt_1 \quad (4.1)$$

Theorem 4
We have the relation

$$J^s(f)(x) = \frac{1}{\Gamma(s)} \int_0^x (x-y)^{s-1} f(y)dy \quad (4.2)$$

$$= \frac{x^s}{\Gamma(s)} \int_0^1 (1-u)^{s-1} f(ux)du \quad (4.3)$$

in which $\Gamma(s)$ is the Euler gamma function for positive integer $s$ [9],[10].

Proof
Let us demonstrate the first relation (4.2). It is evident for $s = 1$. It is true for $s = 2$.

$$J^2(f)(x) = \frac{1}{\Gamma(2)} \int_0^x (x-y) f(y) dy$$

$$= \int_0^x (x-y) f(y) dy \quad (4.4)$$

Let us assume that it is true for $s$ (formula 4.1). We will show that it stands for $s + 1$ too.

$$J^{s+1}(f)(x) = \frac{1}{\Gamma(s)} \int_0^x dt_1 \int_0^{t_1} (t_1-t_2)^{s-1} f(t_2) dt_2 \quad (4.5)$$

We apply the Dirichlet's formula given by Whittaker and Watson [7],[8]

$$\int_a^x dy(x-y)^{\alpha-1} \int_a^y dz(y-z)^{\beta-1} g(y,z)$$
$$= \int_a^x dz \int_z^x dy\, (x-y)^{\alpha-1}(y-z)^{\beta-1} g(y,z) \quad (4.6)$$

for $a = 0, \alpha = 1, \beta = s, y \to t_1, z \to t_2, g(y,z) \to f(t_2)$
Then,

$$J^{s+1}(f)(x) = \frac{1}{\Gamma(s)} \int_0^x dt_2 \int_{t_2}^x dt_1\, (t_1-t_2)^{s-1} f(t_2) \quad (4.7)$$

We integrate over the variable $t_1$

$$J^{s+1}(f)(x) = \frac{1}{\Gamma(s)} \int_0^x dt_2\, (x-t_2)^{s-1} \frac{1}{s} f(t_2)$$

$$= \frac{1}{\Gamma(s+1)} \int_0^x dt_2\, (x-t_2)^s f(t_2) \quad (4.8)$$

Theorem 5
We have the semi group properties

$$J^{s_1} J^{s_2} = J^{s_1+s_2} = J^{s_2} J^{s_1} \quad (4.9)$$

$$D^{s_1} D^{s_2} = D^{s_1+s_2} = D^{s_2} D^{s_1} \quad (4.10)$$

for any positive integers $s_1$ and $s_2$.
Proof
The second relation is immediate. For the first one, we have

$$J^{s_1+s_2}(f)(x)$$
$$= \frac{1}{\Gamma(s_1+s_2)} \int_0^x (x-y)^{s_1+s_2-1} f(y) \quad (4.11)$$

Let us apply $D^1$ to the two sides

$$D^1 J^{s_1+s_2}(f)(x) = \frac{1}{\Gamma(s_1+s_2)} D^1 \int_0^x (x-y)^{s_1+s_2-1} f(y) dy$$

$$= \frac{s_1+s_2-1}{\Gamma(s_1+s_2)} \int_0^x (x-y)^{s_1-2+s_2} f(y) dy$$

$$= J^{s_1-1+s_2}(f)(x) \quad (4.12)$$

We apply successively $D^2, D^3, \ldots D^{s_1}$

$$D^{s_1} J^{s_1+s_2}(f)(x) = \frac{1}{\Gamma(s_2)} D^1 \int_0^x (x-y)^{s_2-1} f(y)$$

$$= J^{s_2}(f)(x)$$

$$D^{s_1} J^{s_1+s_2} = J^{s_2} \quad (4.13)$$

We apply $J^{s_1}$ on the left hand side of the two members

$$J^{s_1} D^{s_1} J^{s_1+s_2} = J^{s_1} J^{s_2} \quad (4.14)$$

$$J^{s_1+s_2} = J^{s_1} J^{s_2} \quad (4.15)$$

This relations is symmetric in $s_1$ and $s_2$:

$$J^{s_1}J^{s_2} = J^{s_1+s_2} = J^{s_2}J^{s_1} \qquad (4.16)$$

## 5. Extension 2. Definitions of S-Order Integral Operator $J^s$ and S-Order Derivative Operator $D^s$ for Any Negative Integer $s$.

We take the advantage that we have demonstrated that the inverse noted $J^{-1}$ of $J^1$ is $D^1$ and the inverse noted $D^{-1}$ of $D^1$ is $J^1$

Then, we define $J^s$ and $D^s$ for any negative integer $s$ by

$$J^s = D^{-s} \qquad D^s = J^{-s} \qquad (5.1)$$

Of course $D^0 = 1_E \quad J^0 = 1_E$. Then for any positive integers $k$ and $s$

$$D^k = J^{-k} \qquad J^k = D^{-k} \qquad (5.2)$$

$$D^k J^s = J^{-k+s} \quad for \ k < s \qquad (5.3)$$

$$D^k J^s = D^{k-p} J^{s-p} \quad for \ any \ p < s \qquad (5.4)$$

## 6. Extension 3. Definitions of S-Order Integral Operator $J^s$ for Any Positive Real $s$ (Fractional, Non Fractional)

We define $J^s$ as

$$J^s(f)(x) = \frac{1}{\Gamma(s)} \int_0^x (x-y)^{s-1} f(y) dy \qquad (6.1)$$

$$= \frac{x^s}{\Gamma(s)} \int_0^1 (1-u)^{s-1} f(ux) du \qquad (6.2)$$

for any positive $s$ and $\Gamma(s)$ is the extension of Euler's gamma function for positive real $s$.

The second relation is obtained from the first one by the change of variable $u = \frac{y}{x}$

**Example 1**
Let us take $f(x) = x^k$ for any positive $k$

$$J^s(x^k) = \frac{x^s}{\Gamma(s)} \int_0^1 (1-u)^{s-1} u^k x^k du$$

$$= \frac{x^{s+k}}{\Gamma(s)} \int_0^1 (1-u)^{s-1} u^k du$$

$$= \frac{x^{s+k}}{\Gamma(s)} B(k+1)(s) \qquad (6.3)$$

in which, we have the extension of Euler's beta function defined by

$$B(p)(q) = \int_0^1 t^{p-1} (1-t)^{q-1} dt \qquad (6.4)$$

for any positive $p$ and $q$ with the property

$$B(p)(q) = \frac{\Gamma(p)\Gamma(q)}{\Gamma(p+q)} \qquad (6.5)$$

Then

$$J^s(x^k) = x^{s+k} \frac{\Gamma(k+1)}{\Gamma(k+s+1)} = x^{s+k} \frac{k!}{(k+s)!} \qquad (6.6)$$

**Example 2**
Let us suppose $s = \frac{1}{2}$ and $k = \frac{1}{2}$

$$J^{\frac{1}{2}}(x^{\frac{1}{2}}) = x^{\frac{1}{2}+\frac{1}{2}} \frac{\Gamma\left(\frac{1}{2}+1\right)}{\Gamma\left(\frac{1}{2}+\frac{1}{2}+1\right)} = \frac{\sqrt{\pi}}{2} x \qquad (6.7)$$

Let us calculate now

$$J^{\frac{1}{2}}(x) = x^{\frac{1}{2}+1} \frac{\Gamma(1+1)}{\Gamma(1+\frac{1}{2}+1)} = \frac{4}{3\sqrt{\pi}} x^{\frac{3}{2}} \qquad (6.8)$$

$$J^{\frac{1}{2}} J^{\frac{1}{2}}(x^{\frac{1}{2}}) = J^{\frac{1}{2}}(\frac{\sqrt{\pi}}{2} x) = \frac{2}{3} x^{\frac{3}{2}} \qquad (6.9)$$

Let us calculate directly

$$J^1(x^{\frac{1}{2}}) = \int_0^x t^{\frac{1}{2}} dt = \frac{2}{3} x^{\frac{3}{2}} \qquad (6.10)$$

So we have from (6.9) and (6.10)

$$J^{\frac{1}{2}} J^{\frac{1}{2}}(x^{\frac{1}{2}}) = J^1(x^{\frac{1}{2}}) \qquad (6.11)$$

**Theorem 6**

$$J^a J^b(x^k) = J^{a+b}(x^k) = J^b J^a(x^k)$$

$$= \frac{\Gamma(k+1)}{\Gamma(a+b+k+1)} x^{a+b+k} \qquad (6.12)$$

for any positive numbers $a, b, k$.
**Proof**

$$J^a J^b(x^k) = J^a(x^{b+k} \frac{\Gamma(k+1)}{\Gamma(k+b+1)})$$

$$= x^{a+b+k} \frac{\Gamma(k+1)}{\Gamma(k+b+1)} \frac{\Gamma(b+k+1)}{\Gamma(k+b+a+1)}$$

$$= \frac{\Gamma(k+1)}{\Gamma(k+b+a+1)} x^{a+b+k} = J^{a+b}(x^k) \qquad (6.13)$$

Then

$$J^a J^b = J^{a+b} = J^b J^a \qquad (6.14)$$

applied to $f(x) = x^k$. This relation is the particular case of general case for $f(x) = x^k$

**Example 3**

Let us suppose $a = \frac{3}{2}$, $k = 2$. Then let us calculate $J^{\frac{3}{2}}(x^2)$. There are many ways to perform this calculation. The simplest one is to split $J^{\frac{3}{2}}$ into $J^{\frac{1}{2}+1}$

$$J^{\frac{3}{2}}(x^2) = J^{\frac{1}{2}}J^1(x^2) = J^{\frac{1}{2}}\int_0^x t^2 dt = \frac{1}{3}J^{\frac{1}{2}}(x^3)$$

$$= \frac{3!}{\frac{7}{2} \times \frac{5}{2} \times \frac{3}{2} \times \frac{1}{2}\Gamma(\frac{1}{2})}\frac{x^{\frac{7}{2}}}{3} = \frac{32}{105\sqrt{\pi}}x^{\frac{7}{2}} \quad (6.15)$$

The direct application of the formula (6.1) gives the result without doing again any integration

$$J^{\frac{3}{2}}(x^2) = \frac{\Gamma(2+1)}{\Gamma(2+\frac{3}{2}+1)}x^{\frac{7}{2}} = \frac{\Gamma(3)}{\Gamma(\frac{9}{2})}x^{\frac{7}{2}} = \frac{32}{105\sqrt{\pi}}x^{\frac{7}{2}} \quad (6.16)$$

Remarkable relations derived from the relation (6.6)

We utilize the expression of $J^s(x^k)$ given in the relation (6.6)

$$J^s(x^k) = x^{s+k}\frac{\Gamma(k+1)}{\Gamma(k+s+1)} = \frac{k!}{k+s!}x^{s+k} \quad (6.17)$$

for any positive $s$ and $k$. We exchange $s$ and $k$

$$J^k(x^s) = x^{s+k}\frac{\Gamma(s+1)}{\Gamma(k+s+1)} = x^{s+k}\frac{s!}{k+s!} \quad (6.18)$$

The ratio

$$\frac{J^s(x^k)}{J^k(x^s)} = \frac{\Gamma(k+1)}{\Gamma(s+1)} = \frac{k!}{s!} \quad (6.19)$$

is independent on $x$. In particular, let us take $s$ and $k$ equal to the transcendental numbers $\pi$ and $e$. Then,

$$J^\pi(x^e) = \frac{\Gamma(e+1)}{\Gamma(e+\pi+1)}x^{\pi+e} = \frac{e!}{(e+\pi)!}x^{\pi+e} \quad (6.20)$$

$$J^e(x^\pi) = \frac{\Gamma(e+1)}{\Gamma(e+\pi+1)}x^{e+\pi} = \frac{\pi!}{(e+\pi)!}x^{\pi+e} \quad (6.21)$$

Then, the ratio

$$\frac{J^\pi(x^e)}{J^e(x^\pi)} = \frac{\Gamma(e+1)}{\Gamma(\pi+1)} = \frac{e!}{\pi!} \quad (6.22)$$

$$= 0.592\,761\,747\,048\,502\,880\,285\,354\,552\,447\,52$$

is independent on $x$. $e!$ and $\pi!$ are the extensions of factorial function for $e$ and $\pi$. This result is remarkable. The different numerical values of $e!$, $(e+\pi)!$, $\pi!$, $\frac{\pi!}{e!}$ and $\frac{e!}{\pi!}$ are given in Appendix 3.

## 7. Definition of $D^s$ for Any Positive $s$ (Fractional, Non Fractional)

We have shown that

$$D^1 J^s(f)(x) = J^{s-1}(f)(x) \quad (7.1)$$

$$D^k J^s(f)(x) = J^{s-k}(f)(x) \quad (7.2)$$

for any positive integer $k$ and any positive integer $s$,. We use this relation to DEFINE $D^s$ for any positive $s$. We CHOOSE an integer $k > s$ such

$$D^s(f)(x) = D^k J^{k-s}(f)(x) \quad (7.3)$$

**Theorem 7**

In the case of $f(x) = x^p$ for any positive $p$ and any positive $s$, the expression

$$D^s(x^p) = D^k J^{k-s}(x^p) = \frac{\Gamma(p+1)}{\Gamma(p-s+1)}x^{p-s}$$

$$= \frac{p!}{(p-s)!}x^{p-s} \quad (7.4)$$

is independent on the choice of $k$. However the intermediate calculations concerning the gamma function require $k$ to be choosen so that $(k-s) > 0$.

Proof.

$$D^k J^{k-s}(x^p) = \frac{D^k}{\Gamma(k-s)}\int_0^x (x-t)^{k-s-1} t^p dt \quad (7.5)$$

We introduce the variable $u = \frac{t}{x}$ in the integral and we remark that

$$\int_0^x (1-u)^{k-s-1} u^p du = B(k-s)(p+1) \quad (7.6)$$

in which $B(k-s)(p+1)$ is Euler's beta function

$$B(k-s)(p+1) = \frac{\Gamma(k-s)\Gamma(p+1)}{\Gamma(k-s+p+1)} \quad (7.7)$$

$$D^k J^{k-s}(x^p) = D^k \frac{1}{\Gamma(k-s)}x^{k-s+p}B(k-s)(p+1)$$

$$= \frac{\Gamma(p+1)}{\Gamma(p-s+1)}x^{p-s} = \frac{p!}{(p-s)!}x^{p-s} \quad (7.8)$$

Remarks

The relation (7.8) is remarkable for many reasons. $D^k J^{k-s}(x^p)$ is independent on the choice of $k$, In practical direct calculation, the result is easily obtained by taking $k$ as small as possible, for instance $k = 1$ or 2.

It is worth noting that even the value of $k$ is arbitrary, we must keep in mind that $k$ is to be choosen in order that all the expressions of intermediate gamma function do exist. It is important to note that the condition $f(0) = 0$ must be fulfilled too.

Example 7.1

Let us calculate $D^{\frac{1}{2}}(x)$ and $k = 1$ in the definition

$$D^{\frac{1}{2}}(x) = D^1 J^{\frac{1}{2}}(x) = D^1(x^{\frac{3}{2}}) \frac{1}{\Gamma(\frac{1}{2})} B(\frac{1}{2}, 2) = \frac{2}{\sqrt{\pi}} x^{\frac{1}{2}} \quad (7.9)$$

Let us choose now $k = 2$

$$D^{\frac{1}{2}}(x) = D^2 J^{\frac{3}{2}}(x) = D^2(x^{\frac{5}{2}}) \frac{1}{\Gamma(\frac{3}{2})} B\left(\frac{3}{2}, 2\right) = \frac{2}{\sqrt{\pi}} x^{\frac{1}{2}} \quad (7.10)$$

Let us apply the formula (7.4)

$$D^{\frac{1}{2}}(x) = \frac{\Gamma(1+1)}{\Gamma(1 - \frac{1}{2} + 1)} x^{1-\frac{1}{2}} = \frac{\Gamma(2)}{\Gamma(\frac{3}{2})} x^{\frac{1}{2}} = \frac{2}{\sqrt{\pi}} x^{\frac{1}{2}} \quad (7.11)$$

Of course, we obtain the same result for the 3 calculations

$$D^{\frac{1}{2}} D^{\frac{1}{2}}(x) = D^1(x) = 1 \quad (7.12)$$

Example 7.2
Let us apply $D^{\frac{1}{2}}$ to the function $f(x) = x^{-\frac{1}{2}}$. The latter one is infinite and is not equal to zero for $x = 0$. Then, we are expecting to meet a difficulty.

$$D^{\frac{1}{2}}(x^{-\frac{1}{2}}) = \frac{\Gamma\left(-\frac{1}{2} + 1\right)}{\Gamma(-1 + 1)} x^{-\frac{1}{2} - \frac{1}{2}} = \frac{\Gamma\left(\frac{1}{2}\right)}{\Gamma(0)} x^{-1} = 0 \quad (7.14)$$

because $\Gamma(0)$ is infinite. Then

$$D^{\frac{1}{2}} D^{\frac{1}{2}}(x^{-\frac{1}{2}}) = D^{\frac{1}{2}}(0) = 0 \quad (7.15)$$

But

$$D^1(x^{-\frac{1}{2}}) = -\frac{1}{2} x^{-\frac{3}{2}} \neq 0 \quad (7.16)$$

Then

$$D^{\frac{1}{2}} D^{\frac{1}{2}}(x^{-\frac{1}{2}}) \neq D^1(x^{-\frac{1}{2}}) \quad (7.17)$$

Theorem 8
If $f(x) = x^p$ for any positive $p$ and any positive $a$ and $b$, then
$$D^a D^b(x^p) = D^b D^a(x^p) = D^{a+b}(x^p) \quad (7.18)$$
Proof
Let us utilize the result of the Theorem 7

$$D^a D^b(x^p) = D^a \frac{\Gamma(p+1)}{\Gamma(p - b + 1)} x^{p-b} \quad p > b$$

$$= \frac{\Gamma(p+1)}{\Gamma(p - b + 1)} D^a(x^{p-b}) \quad a > p - b$$

$$= \frac{\Gamma(p+1)}{\Gamma(p - b - a + 1)} x^{p-b-a}$$

$$= D^{a+b}(x^p) = D^b D^a(x^p) \quad (7.19)$$

for $p - a - b > 0$
Remarks
1) If we take $s = e$ and $p = \pi$, then

$$D^e(x^\pi) = \frac{\Gamma(\pi + 1)}{\Gamma(\pi - e + 1)} x^{\pi - e} = \frac{\pi!}{(\pi - e)!} x^{\pi - e} \quad (7.20a)$$

$$D^\pi(x^e) = \frac{\Gamma(e + 1)}{\Gamma(e - \pi + 1)} x^{e - \pi} = \frac{e!}{(e - \pi)!} x^{e - \pi} \quad (7.20b)$$

$$D^e(x^\pi) \cdot D^\pi(x^e) = \frac{\pi!}{(\pi - e)!} \frac{e!}{(e - \pi)!}$$

$$= 22.364\,994\,517\,058\,857\,454\,906\,921\,720\,114$$

This expression may be derived from the relation (6.20).
2) The formulae (7.20a) and (7.20b) may be deduced from the relation (6.20) and (6.21) by changing $\pi$ into $-\pi$, $e$ into $-e$ and taking account of $J^{-\pi} = D^\pi$ and $J^{-e} = D^e$
The numerical values of $\frac{\pi!}{(\pi - e)!}$ and $\frac{e!}{(e - \pi)!}$ are given in the Appendix 3.
3) If in the relation

$$J^s(x^p) = \frac{\Gamma(p + 1)}{\Gamma(s + p + 1)} x^{s+p} \quad (7.21)$$

we change the positive number $s$ into $-s$, we obtain exactly the expression (7.4) of $D^s(x^p)$. This justifies our remark $D^s = J^{-s}$ or $D^{-s} = J^s$ for any positive s (Theorem 2).
It is worth pointing out once more that the condition $f(x) = 0$ for $x = 0$ is essential.
Theorem 9
If $f$ belong to the set $E$ defined in the introduction $(f(0) = 0)$, then for any positive $a$ and $b$, we have

$$D^a D^b(f)(x) = D^{a+b}(f)(x) = D^b D^a(f)(x) \quad (7.22)$$

Proof.
Let us choose the integer $k$ such as $k > a + b$

$$D^a D^b(f)(x) = D^a D^k J^{k-b}(f)(x) \quad k > b \text{ and } k \in \mathbb{N}$$

$$= D^a D^k \frac{1}{\Gamma(k - b)} \int_0^x (x - y)^{k-b-1} f(y) dy \quad (7.23)$$

We utilize Fubini's theorem noting that the contributions of the $(x - y)^{k-b-1} f(y)$ are all zero for $y = x$ and $y = 0$ because $f(0) = 0$.

$$D^a D^b(f)(x) = \frac{D^a}{\Gamma(k - b)} \int_0^x (x - y)^{k-b-1} f(t) dt$$

$$= \frac{1}{\Gamma(k - b - a)} \int_0^x (x - y)^{k-b-a-1} f(t) dt$$

$$= D^k J^{k-b-a}(f)(x)$$

$$= D^{a+b}(f)(x) = D^b D^a(f)(x) \quad (7.24)$$

Then

$$D^a D^b = D^{a+b} = D^b D^a \quad (7.25)$$

if it is applied to the set $E$ of function $f(x) = 0$ for $x = 0$. The relation is not valid if the latter condition is not fulfilled. Let us give an illustrative example. Let us take $f(x) = constant = 1$ for every $x$ ($f(0) = 1 \neq 0$).

$$D^{\frac{1}{2}}(1) = D^1 J^{1-\frac{1}{2}}(f)(x) = D^1 J^{\frac{1}{2}}(f)(x)$$
$$= \frac{1}{\sqrt{\pi}} x^{-\frac{1}{2}} \quad (7.26)$$

We find that the half derivative of a constant function is different of zero and depends on $x$ which is an absurd result. If $f(x) = constant = 0$ instead of 1 for every $x$, then

$$D^{\frac{1}{2}}(0) = 0 \qquad D^{\frac{1}{2}}D^{\frac{1}{2}}(0) = D^1(0) = 0$$

which is not a nonsense because $f(x) = 0$ for any $x$, then for $x = 0$ too.

## 8. Conclusion

We have proposed three definitions of $s$-order derivative $D^s$ and $k$-order integral $J^k$. All the three are equivalent for positive integer $s$ and positive integer $k$.

It is shown that $D^s = J^{-s}$ and $J^k = D^{-k}$ for any positive integer $s$ and $k$. The extension of the definition of $D^s$ and $J^k$ for any positive number $s$ and $k$ is given too.

Applying the operators $J^s$ and $D^k$ to the function $x^p$ for any positive $p$, we obtain the following remarkable results

$$J^s(x^p) = x^{s+p} \frac{\Gamma(p+1)}{\Gamma(p+s+1)} = \frac{p!}{(k+p)!} x^{s+p} \quad (8.1)$$

for any positive $s$ and $p$ and

$$D^k(x^p) = D^n J^{n-k}(x^p) = \frac{\Gamma(p+1)}{\Gamma(p-k+1)} x^{p-k} \quad (8.2)$$

for any positive numbers $k, p$ and $n$; $n > k$, . The result is independent on the choice of $n$.

If $s, p$ and $k$ are equal to transcendental numbers $\pi$ and $e$, we obtain

$$J^\pi(x^e) = \frac{\Gamma(e+1)}{\Gamma(e+\pi+1)} x^{\pi+e} = \frac{e!}{(e+\pi)!} x^{\pi+e} \quad (8.3)$$

$$J^e(x^\pi) = \frac{\Gamma(e+1)}{\Gamma(e+\pi+1)} x^{e+\pi} = \frac{\pi!}{(e+\pi)!} x^{\pi+e} \quad (8.4)$$

$$\frac{J^\pi(x^e)}{J^e(x^\pi)} = \frac{\Gamma(e+1)}{\Gamma(\pi+1)} = \frac{e!}{\pi!} \quad (8.5)$$

The ratio is independent on $x$.

$$D^e(x^\pi) = \frac{\Gamma(\pi+1)}{\Gamma(\pi-e+1)} x^{\pi-e} = \frac{\pi!}{(\pi-e)!} x^{\pi-e} \quad (8.6)$$

$$D^\pi(x^e) = \frac{\Gamma(e+1)}{\Gamma(e-\pi+1)} x^{e-\pi} = \frac{e!}{(e-\pi)!} x^{e-\pi} \quad (8.7)$$

$$D^e(x^\pi) \cdot D^\pi(x^e) = \frac{e! \, \pi!}{(\pi-e)! \, (e-\pi)!} \quad (8.8)$$

The product is independent on x.

$\Gamma(p)$ is the extension of the Euler's gamma function for any number $p$ (not necessary integer number) and $p! = \Gamma(p + 1)$ is the extension of factorial function for any positive number $p$. The numerical values concerning $\pi$ and $e$, $\pi!$ and $e!$ correct to 31 decimals are given in the Appendix 3.

As a final conclusion, we would like to say that the approach by means of integral operators and derivative operators derived from integral operators is better and more general than the two first ones that we have given in our work [5] because the properties of the set $E$ of functions $f$ we are looking for (for instance causal function of the real variable $x$) are assumed to be known from the beginning. The limit of the applicability of the method is then clear.

We give the study of the case of complex order integrals and derivatives in our work [11].

## Appendix

### *Appendix 1: The Euler's gamma and beta functions*

The gamma function is defined by the integral

$$\Gamma(p) = \int_0^\infty t^{p-1} e^{-t} dt \quad (A1.1)$$

for any real positive $p$. It is easy to verify that $\Gamma(p)$ has the property

$$\Gamma(p+1) = p\Gamma(p) \quad (A1.2)$$

For positive integer $p$

$$\Gamma(p+1) = p! \quad (factorial \, p) \quad (A1.3)$$

for any positive $p$ (not necessary integer); we may define then a generalization of the factorial

$$p! = \Gamma(p+1) \quad (A1.5)$$

So we know the meaning of the expressions $\frac{p}{q}!$, $e!$, $\pi!$ for the fraction $\frac{p}{q}$, the transcendental numbers $e$ and $\pi$

The beta function is defined by

$$B(x)(y) = \int_0^1 t^{x-1} (1-t)^{y-1} dt \quad (A1.6)$$

for any positive real numbers $x$ and $y$. $B(x)(y)$ is symmetric in $x$ and $y$

$$B(x)(y) = \frac{\Gamma(x)\Gamma(y)}{\Gamma(x+y)} \quad (A1.7)$$

We have shown

$$J^s(f)(x) = \frac{1}{\Gamma(s)} \int_0^x (t-x)^{s-1} f(t) dt \quad (A1.8)$$

If we have

$$f(x) = x^k \quad (A1.9)$$

for any positive $k$, then

$$J^s(x^k) = \frac{1}{\Gamma(s)} \int_0^x (t-x)^{s-1} t^k dt \quad (A1.10)$$

We introduce the variable $u = t/x$, then

$$J^s(x^k) = \frac{x^{s+k}}{\Gamma(s)} \int_0^x u^k (1-u)^{s-1} du = \frac{x^{s+k}}{\Gamma(s)} B(k+1)(s)$$

$$= \frac{x^{s+k}}{\Gamma(s)} \frac{\Gamma(k+1)\Gamma(s)}{\Gamma(s+k+1)} = \frac{k!}{(s+k)!} x^{s+k} \quad (A1.11)$$

**Appendix 2: About the inverse of an operator [5]**

Theorem
Any operator has at least a right hand side inverse $C$
Proof.
Let $y$ be an element of the value domain $Val(A)$ of an operator $A$. By definition of $Val(A)$ there is at least an element $x_0$ belonging to the definition domain $Def(A)$ of the operator $A$ such as

$$A(x_0) = y \quad (A2.1)$$

for any element $y$ of $Val(A)$, we may choose an $x$. Let us designate it by $x_0$ and define the operator $C$ such as

$$C(y) = x_0 \quad (A2.2)$$

Then $A(C)(y) = A(x_0) = y$ for any $y \in Val(A)$

$$AC = 1_{Val(A)} \quad (A2.3)$$

in which $1_{Val(A)}$ is the identity operator on $Val(A)$. $C$ is a right hand side inverse of $A$. It depends on a choice.
Examples
The inverse function in the classical meaning is in fact a right hand side inverse. As an example, Arcsin is the right hand side inverse function of the sinus function

$$sin(Arcsin) = 1_{Val(sin)} = 1_{[-1,1]}$$

in which $[-1, 1]$ is the segment of the line between $-1$ to $1$ for any $a \in [-1,1]$.

$$sin(Arcsin)(a) = sin\beta = a$$

in which $\beta = Arcsin(a)$ is the principal determination belonging to the angle $[-\frac{\pi}{2}, \frac{\pi}{2}]$.
Let us now look for the meaning of $Arcsin(sin)$ for any real number (any angle) $x \in \mathbb{R} = Def(sin)$

$$sinx = \alpha$$

$$Arcsin(\alpha) = x_0 \in [-\frac{\pi}{2}, \frac{\pi}{2}]$$

$$Arcsin(sin)(x) = Arcsin(sinx)$$

$$Arcsin(\alpha) = x_0$$

Then we have

$$Arcsin(sin)(x) \neq x_0 \text{ or } Arcsin(sin) \neq 1_{Defsin}$$

$$Val(Arcsin(sin)) = [-\frac{\pi}{2}, \frac{\pi}{2}]$$

$$Def(Arcsin(sin)) = Def(sin) = \mathbb{R}$$

Let us point out that for the inverse functions we have the problem of choice. The inverse function implies a choice for the determination to be choosen for the functions having many determinations (functions like $\sqrt[n]{\phantom{x}}$ , arcsin, arctan, etc…) It is worth noting that $Arcsin$ is an operator, but $arcsin$ is not an operator because there are many (infinite) values of $arcsin(x)$ for a given $x$.
The definition of the derivative operator $D^s$ obtained from the definition of the integral operator $J^s$

$$D^s(f)(x) = D^k J^{k-s}(f)(x)$$

introduces the choice on the positive number $k$ ($k - s > 0$). We have shown that the final result is in fact independent on **k**. We have the following relation

$$J^1(f)(x) = \int_a^x f(t) dt = F(x) - F(a)$$

in which $F(x)$ is a primitive of $f(x)$

$$D^1 J^1(f)(x) = F'(x) = f(x) \text{ or } D^1 J^1 = 1_E$$

$$J^1 D^1(f)(x) = \int_0^x f'(t) dt = f(t) - f(0) \neq f(x)$$

$$J^1 D^1 \neq 1_E$$

$$J^1 D^1 = 1_E \text{ if and only if } f(0) = 0$$

**Appendix 3: Numerical values for the transcendental numbers $\pi$ and $e$**

$e = 2.718\ 281\ 828\ 459\ 045\ 235\ 360\ 287\ 471\ 352\ 7$

$e! = 4.260\ 820\ 476\ 357\ 003\ 381\ 700\ 121\ 224\ 647\ 7$

$\pi = 3.141\ 592\ 653\ 589\ 793\ 238\ 462\ 643\ 383\ 279\ 5$

$\pi! = 7.188\ 082\ 728\ 976\ 032\ 702\ 082\ 194\ 345\ 124\ 8$

$\frac{e!}{\pi!} = 0.592\ 761\ 747\ 048\ 502\ 880\ 285\ 354\ 552\ 437\ 32$

$\frac{\pi!}{e!} = 1.687\ 018\ 443\ 715\ 759\ 455\ 687\ 799\ 928\ 242\ 6$

$(\pi - e) = 0.423\ 310\ 825\ 130\ 748\ 003\ 102\ 355\ 911\ 926\ 84$

$(\pi - e)! = 0.886\ 240\ 147\ 692\ 794\ 559\ 514\ 959\ 131\ 20817$

$(e - \pi) = 0.423\ 310\ 825\ 130\ 748\ 003\ 102\ 355\ 911\ 926\ 84$

$(e - \pi)! = 1.545\ 204\ 936\ 151\ 901\ 746\ 654\ 139\ 877\ 849\ 1$

$\dfrac{\pi!}{(\pi - e)!} = 8.110\ 761\ 792\ 601\ 279\ 051\ 112\ 802\ 855\ 137\ 1$

$\dfrac{e!}{(e - \pi)!} = 2.757\ 446\ 845\ 185\ 422\ 310\ 617\ 384\ 631\ 128\ 6$

$\dfrac{\pi!}{(\pi - e)!}\dfrac{e!}{(e - \pi)!} = 22.364\ 994\ 517\ 058\ 857\ 454\ 906\ 921\ 720\ 114$

$(\pi + e) = 5.859\ 874\ 482\ 048\ 838\ 473\ 822\ 930\ 854\ 632\ 2$

$(\pi + e)! = 554.654\ 105\ 737\ 269\ 399\ 798\ 013\ 158\ 641\ 18$

$\dfrac{\pi!}{(\pi + e)!} = 0.012\ 959\ 577\ 247\ 555\ 632826589\ 943903911$

## References

[1] S. Miller, Kenneth, "An introduction to the fractional calculus and the fractional differential equations", Bertram Ross (Editor). Publisher: John Wiley and Sons 1st edition, 1993, ISBN 0-471-58884-9.

[2] B. Oldham Keith and Spanier Jerome, "The fractional calculus. Theory and Application of differentiation and integration to arbitrary order" (Mathematics in Science and engineering). Publisher: Academic Press, Nov 1974, ISBN 0-12-525550-0.

[3] Zavada, "Operator of fractional derivative in the complex plane", Institute of Physics, Academy of Sciences of Czech Republic, 1997.

[4] F. Dubois, A. C Galucio, N. Point, "Introduction à la dérivation fractionaire. Théorie et Applications", http://www.mathu-psud.fr/~fdubois/travaux/evolution/ananelly-07/techinge-2010/dgp-fractionnaire-mars2010.pdf, 29 Mars 2010.

[5] Raoelina Andriambololona, Tokiniaina Ranaivoson, Rakotoson Hanitriarivo, Roland Raboanary, , "Two definitions of fractional derivatives of power functions". Institut National des Sciences et Techniques Nucleaires (INSTN-Madagascar), 2012, arXiv:1204.1493.

[6] Raoelina Andriambololona, "Algèbre linéaire et multilinéaire." Applications. 3 tomes. Collection LIRA-INSTN Madagascar, Antananarivo, Madagascar, 1986 Tome 1, pp 2-59.

[7] E.T. Whittaker, and G.N. Watson, "A course of modern analysis", Cambridge University Press, Cambridge, 1965.

[8] R. Herrmann, "Fractional Calculus. An introduction for Physicist", World Scientific Publishing, Singapore, 2011.

[9] E. Artin, "The Gamma Function", Holt, Rinehart and Winston, New York, 1964.

[10] S.C. Krantz, "The Gamma and beta functions" § 13.1 in handbook of complex analysis, Birkhauser, Boston, MA, 1999, pp.155-158.

[11] Raoelina Andriambololona, Tokiniaina Ranaivoson, Rakotoson Hanitriarivo, "Definition of complex order integrals and derivatives using operator approach", INSTN preprint 120829, arXiv 1409.400. Published in IJLRST, Vol.1,Issue 4:Page No.317-323, November-December(2012), ISSN (online):2278-5299, http://www.mnkjournals.com/ijlrts.htm.